\documentclass[11pt]{amsart}

\usepackage{a4wide}
\usepackage{amsmath,amsfonts,amsthm,amssymb}

\def\d{{\mathrm{d}}}

\def\D{{\mathcal{D}}}
\def\O{{\mathcal{O}}}
\def\bydef{:=}

\newcommand{\rd}{\mathrm{d}}
\newcommand{\re}{\mathrm{e}}

\newcommand{\bm}{\mathbf{\mathcal{BM}^+}([0,1])}
\newcommand{\bmi}{\mathbf{\mathcal{BM}^+}((0,1))}

\newcommand{\bmg}{\mathbf{\mathcal{BM}}([0,1])}
\newcommand{\bmgi}{\mathbf{\mathcal{BM}}((0,1))}
\newcommand{\lbm}{L^{\infty}\left([0,\infty);\bm\right)}
\newcommand{\lbmg}{L^{\infty}\left([0,\infty);\bmg\right)}

\newcommand{\cicic}{C^{\infty}\left(\rb^+;C^{\infty}([0,1])\right)}
\newcommand{\len}{L^{\infty}\left([0,\infty);\mathcal{E}'\right)}

\newcommand{\eps}{\epsilon}
\newcommand{\supp}{\mathrm{supp}}
\newcommand{\ssupp}{\textrm{sing supp}}

\newcommand{\rb}{\mathbb{R}}
\newcommand{\bt}{\widetilde{\mathcal{T}}}
\newcommand{\bti}[2]{\widetilde{\mathcal{T}}_{#1,#2}}
\newcommand{\E}{\mathcal{E'}}
\newcommand{\hvp}{\hat{\varphi}}

\newtheorem {thm} {Theorem}
\newtheorem {lem} {Lemma}
\newtheorem {rmk}  {Remark}

\newtheorem {defi} {Definition}

\newtheorem{prop}{Proposition}

\newcommand{\eliminado}[1]{{}}

\begin{document}

\title[A non-standard evolution problem]{
A non-standard evolution problem arising in population genetics
}

\author{Fabio A. C. C. Chalub}
\address{Departamento de Matem\'atica and Centro de Matem\'atica
e Aplica\c c\~oes, Universidade Nova de Lisboa, 
Quinta da Torre, 2829-516, Caparica, Portugal.}

\email{chalub@fct.unl.pt}
\author{Max O. Souza}\
\address{Departamento de Matem\'atica Aplicada, Universidade Federal
Fluminense, R. M\'ario Santos Braga, s/n, 22240-920, Niter\'oi, RJ, Brasil.}

\email{msouza@mat.uff.br}

\thanks{%
The authors want to thank Peter Markowich for helpful comments, and an anonymous referee for pointing out a mistake in an earlier version of Proposition 4 and for general comments that improved the presentation.
FACCC is partially supported by FCT/Portugal, grants POCI/MAT/57546/2004, 
PTDC/MAT/68615/2006 and PTDC/MAT/66426/2006. 
MOS is partially supported by FAPERJ/Brazil grants  170.382/2006 and 110.174/2009. MOS also thanks the support and hospitality of FCT/UNL and Complexo Interdisciplinar/UL. Part of this work  has been done during the \textit{Special Semester on Quantitative Biology Analyzed by Mathematics}, organized by RICAM, Austrian Academy of Sciences. 
}



\begin{abstract}
We study the evolution of the probability density of an asexual, one locus population under 
natural selection and random evolution. This evolution is governed by a Fokker-Planck equation with degenerate coefficients on the boundaries, supplemented by a pair of conservation laws. It is readily shown that no classical or standard weak solution  definition yields solvability of the problem. We provide an appropriate definition of weak solution for the problem, for which we show existence and uniqueness. The solution displays a very distinctive structure and, for large time, we show convergence to a unique stationary solution that turns out to be a singular measure supported at the endpoints. An exponential rate of convergence to this steady state is also proved.
\end{abstract}

\maketitle

\section{Introduction}

A classical problem in population genetics is to study the evolution of a mutant gene. 
A standard approach to this problem is to consider a finite size
population and to define a discrete dynamics  for the evolution of the probability density of such a population. Usually, such models are Markov chains, in which the only absorbing states are the two pure ones. Therefore, one expects, for large time, convergence to one of these two states and, depending on which state is achieved, one says that the mutation has been either fixed or lost.  For large populations, it is natural to ask for a continuous model that approximates this evolution.
In a number of different ways, one arrives at a Fokker-Plank equation that describes either the evolution of the probability density (the so-called forward Kolmogorov equation), or what is sometimes called the transient fixation probability (the backward Kolmogorov equation). From a mathematical point of view, it is interesting to notice that, for the fixation probabilities, it is easy to specify the appropriate non-homogeneous boundary conditions  which, after subtraction of an appropriate multiple of a stationary solution, are recast as Dirichlet conditions. Nevertheless, this does not seem to be the case for the probability evolution. Since it must conserve mass, in many cases a condition of null probability current at the endpoints is used (e.g. \cite{McKaneWaxman07}). For a thorough introduction to the several aspects of mathematical population genetics, we refer the reader  to the monographs by  \cite{Burger2000,Ewens} 

For a class of problems, however, these Fokker-Plank equations turn out to have degenerate coefficients at the boundaries, the classical Kimura equation (cf. \cite{Kimura}) being the archetypal example. For the backward one, this is not a problem since the infinitesimal generator is, very generally, self-adjoint.   For the forward equation, however, the underlying spectral problem is of the limit-point type and, thus, no boundary conditions can be enforced. In particular, one cannot control the flux of the solutions across the boundary of the corresponding domain, and the existence of conservation laws are not to be expected in general. This is an old issue in the study of diffusions, and  it has been tackled  by \cite{Feller1952}, where the so called lateral conditions are derived, in order to ensure that the forward and backward equations are adjoint to each other. With these lateral conditions, however, the forward equation looses its differential character, and this led to a prevalence of the backward equation in the study of diffusions (particularly after \cite{Feller1954}). We shall see below that is possible to  ensure the duality of the backward and forward equations, while maintaining the differential character of the forward equation,  within the framework of weak solutions.

We shall study the forward Kolmogorov equation
\begin{equation}
 \left\{
\begin{array}{lr}
\partial_t p(t,x) = \partial_x^2\left(F(x)p(t,x)\right) - 
\partial_x\left(G(x)p(t,x)\right),&x\in(0,1),\quad t>0\\
p(0,x)=p^{0}(x)
\end{array}
\right.
\label{our:pde}
\end{equation}
with $F$ positive in $(0,1)$, but with simple zeros at the endpoints, and with $G$ vanishing at the endpoints\footnote{More precise statements on the hypothesis made upon $F$ and $G$ are deferred to section~\ref{prelim}.}. Typical examples are $F(x)=G(x)=x(1-x)$ (forward Kimura) and $F(x)=x(1-x)$, $G(x)=x(1-x)(\eta x +\beta)$ (forward Kimura with frequency selection; see \cite{ChalubSouza06a}). 

Equation (\ref{our:pde}) is supplemented by the following conservation laws:
\begin{subequations}
\begin{align}
&\frac{\d}{\d t}\int_0^1p(t,x)\,\d x = 0,\label{cons:prob}\\
&\frac{\d}{\d t}\int_0^1\psi(x)p(t,x)\,\d x = 0,\label{cons:psi}
\end{align}
\label{cons:laws}
\end{subequations}
where $\psi$ satisfies
\begin{equation}
F(x)\psi'' + G(x)\psi' =0,\qquad \psi(0)=0,\quad\text{and}\quad \psi(1)=1.
\label{fixprob}
\end{equation}

\begin{rmk}
In population genetics, the function $\psi$ is referred to as the fixation probability.
 Condition \eqref{cons:prob} is usually stated (or assumed), in the literature of population genetics, but condition \eqref{cons:psi} is not. These conditions have been derived in \cite{ChalubSouza06a}, when obtaining the forward Kimura equation, with frequency selection, as a large population limit of the so called Moran process (cf. \cite{Moran}). See also \cite{TraulsenClaussenHauert} for an alternative approach.
\end{rmk}

Before we proceed, we want to clarify the nature of the conservation laws given by \eqref{cons:laws}. The backward equation and  (formal) adjoint of \eqref{our:pde} is given by
\begin{equation}
 \left\{
\begin{array}{lr}
\partial_t f =F(x)\partial_x^2 f(t,x) + G(x)\partial_x f(t,x)&x\in(0,1),\quad t>0\\
f(0,x)=f^0(x).
\end{array}
\right.
\label{our:adjoint}
\end{equation}
It is readily seen that any stationary solution to \eqref{our:adjoint} is a linear combination of
a constant and  $\psi(x)$. Therefore, the  conservation laws \eqref{cons:laws} are related to the kernel of the infinitesimal generator of \eqref{our:adjoint}. Finally, it should be mentioned that, if (\ref{our:pde}) is a correct approximation of the biological process, then one expects that the probability mass accumulates at the endpoints, as $t$ goes to infinity~\cite{Kimura}.  

The goal of this work is to clarify in what sense a solution to (\ref{our:pde}) that satisfies (\ref{cons:laws}) exists, and  how it behaves for large time. In contradistinction with \cite{Feller1952}, which uses classical function spaces and has to modify equations \eqref{our:pde} and \eqref{our:adjoint} in order to obtain the duality relation, we shall always work with these equations, but in more general, non-normed, distributional spaces. This also differs from recent work in degenerate equations, as for instance: the controlability of degenerate heat equations \cite{Martinez2006}, with solutions in weighted Sobolev spaces; entropy solutions of Fokker-Planck from multilane traffic flow \cite{Dolb1}, where the conditions that might lead to concentration at the end points are explicitly avoided; and from the qualitative  studies by \cite{Dolb2}. See~\cite{Friedman2} for a general discussion of degenerate diffusion equations. We mention also the classical monographs~\cite{CS76,DiBenedetto93}.

 Equation (\ref{our:pde}), with $F(x)=x(1-x)$ and $G(x)=x(1-x)(\eta x +\beta)$ has been studied in reference \cite{ChalubSouza06a}, where a proof of existence and uniqueness in the sense of definition~\ref{our_wsd} is given, under the assumption of interior regularity. An announcement that also includes other results was made in \cite{ChalubSouza06b}. See also \cite{ChalubSouza08b}. More recently, the same problem has been studied through skillful, but formal, calculations in \cite{McKaneWaxman07}, with conditions of null probability current (formally) imposed. Thus, this work can be seen as complementary to the work by \cite{Feller1952} by giving a differential formulation to the forward-backward duality for degenerate diffusions. Also, it can be regarded as an extension of \cite{ChalubSouza06a,ChalubSouza06b}, and as a rigorous proof of the formal calculations in \cite{McKaneWaxman07}.

The main results of the paper can be outlined as follows: let $\bm$ denote the space of (positive) Radon measures on $[0,1]$. Then we have
\begin{thm}[outline]
\label{main_thm}
For a given $p^0\in\bm$, there exists a unique solution $p$ to Equation~(\ref{our:pde}), in a sense to be made precise in definition~\ref{our_wsd}, with $p\in\lbm$ and such that $p$ satisfies the conservations laws (\ref{cons:laws}).  The solution can be written as
\[
p(t,x)=q(t,x)+a(t)\delta_0+b(t)\delta_1,
\]
where $\delta_y$ denotes the singular measure supported 
at $y$, and $q\in\cicic$ is a classical solution to (\ref{our:pde}). We also have that 
$a(t)$ and $b(t)$,
belong to $C([0,\infty))\cap C^{\infty}(\rb^+)$. In particular, we have that 
\[
p\in C^{\infty}(\rb^{+};\bm)\cap C^{\infty}(\rb^{+};C^{\infty}((0,1))\ . 
\]
For large time, we have that $\lim_{t\to\infty}q(t,x)=0$, uniformly, 
and that $a(t)$ and $b(t)$ are monotonically increasing functions such that:
\begin{align*}
&a^\infty\bydef\lim_{t\to\infty} a(t)= \int_0^1(1-\psi(x))p^0(x)\,\rd x\quad\text{and}\\
&b^\infty\bydef\lim_{t\to\infty} b(t)=\int_0^1\psi(x)p^0(x)\,\rd x,
\end{align*}
Moreover, we have that
\[
\lim_{t\to\infty}p(t,\cdot)=a^\infty\delta_0+b^\infty\delta_1,
\] 
with respect to the Radon metric. Finally, the convergence rate is exponential.
\end{thm}

\begin{rmk}
The coefficients of the singular measures, $a(t)$ and $b(t)$ are, respectively, the extinction and the fixation probabilities. Also, notice that the decomposition of $p$ does not follows immediately from the linearity of (\ref{our:pde}). As a matter of fact, neither of the summands are, per se, a solution to (\ref{our:pde}) in the sense of definition~\ref{our_wsd}. Heuristically, as (\ref{our:pde}) is uniformly parabolic in each proper compact set of the unit interval, the parabolic operator erodes  the interior density of the initial measure, which is then transferred into the boundaries and absorbed by the singular measures there. 
\end{rmk}

The outline of the paper is as follows: in section~\ref{prelim}, we present background results for the classical (in a broad sense) solutions to (\ref{our:pde}). In section~\ref{formulation}, we introduce an appropriate definition of a weak solution and show that any solution of this type must satisfy the conservation laws (\ref{cons:laws}). We also show that, with this formulation, (\ref{our:adjoint2}) is indeed the adjoint of (\ref{our:pde}). In section~\ref{the:theorem}, we present the proofs of existence and uniqueness. Section~\ref{largetime} discusses the convergence to the measures supported at the endpoints as time goes to infinity.

\section{Preliminaries}
\label{prelim}

For the convenience of the reader, we present in this section some material
that will be useful in the sequel.

Let $F,G:[0,1]\to\rb$ be smooth, and assume that
\begin{enumerate}
 \item $F$ has single zeros at $x=0$ and at $x=1$, and $F(x)>0$, for $x\in(0,1)$;
\item $G$ has zeros at $x=0$ and $x=1$.
\end{enumerate}
Hadamard's lemma (cf. \cite{BruceGiblin92}) then yields
\[
 F(x)=x(1-x)\Psi(x),\quad \Psi(x)>0\text{ for }x\in[0,1]
\qquad\text{and }
G(x)=x(1-x)\Pi(x)
\]
Let us write, 
\[
 \Xi(x)=\frac{\Pi(x)}{\Psi(x)}.
\]
Then we can rewrite \eqref{our:adjoint} as 
\begin{equation}
 \left\{
\begin{array}{lr}
\partial_t f =x(1-x)\Psi(x)\left[ \partial_x^2 f(t,x) + \Xi(x)\partial_x f(t,x)\right]&x\in(0,1),\quad t>0\\
f(0,x)=f^0(x)
\end{array}
\right.
\label{our:adjoint2}
\end{equation}
The stationary solutions of (\ref{our:adjoint2}) are linear combinations of a constant and
\[
 \psi(x)=c^{-1}\int_0^x\re^{-\int_0^s\Xi(r)\,\rd r}\,\rd s,\quad c=\int_0^1\re^{-\int_0^s\Xi(r)\,\rd r}\,\rd s.
\]

Existence of classical solutions to \eqref{our:pde} can be established by Fourier series, and this is easier done by writing \eqref{our:pde} in selfadjoint form. Let 
\begin{equation}
\re^{\frac{1}{2}\int_0^x\Xi(s)\,\rd s}w=x(1-x)\Psi(x)p.\label{our_change}
\end{equation}
Then (\ref{our:pde}) becomes
\begin{equation}
\partial_tw=x(1-x)\Psi(x)\left\{\partial^2_xw-
\frac{1}{4}\left[2\Xi'+\Xi^2\right]w\right\}.
\label{our:pde:sa}
\end{equation}
\begin{rmk}
Since the  standard maximum principle holds for $C^{1,2}$ solutions of (\ref{our:pde:sa}), we find that, if the initial condition is nonnegative, then $w(t,\cdot)$ is also nonnegative. Moreover, this holds also for $p(t,\cdot)$.
\end{rmk}
 Consider the associated spectral problem:
\begin{equation}
\begin{array}{c}
 -\varphi''+\frac{1}{4}\left[2\Xi'+\Xi^2\right]\varphi=\lambda \theta(x)\varphi,\\
\\
\varphi(0)=\varphi(1)=0, \quad \theta(x)=\frac{1}{\Psi(x)x(1-x)}.
\end{array}
\label{our_eigenop}
\end{equation}
Sturm-Liouville theory for singular problems allows us to conclude that   (\ref{our_eigenop}) is a self-adjoint operator in $L^2\left([0,1],\theta(x)\rd x\right)$, with a complete set of eigenfunctions. In what follows, all $L^2$ spaces will be with respect to $\theta(x)\rd x$, and we shall write $(\cdot,\cdot)$ and $\|\cdot\|_2$ for the corrresponding inner product and norm, respectively. We also  recall, see \cite{CL55,Pryce93} for instance, that
\[
\lim_{j\to\infty}\frac{\lambda_j}{j^2}=K.
\]
An important property of (\ref{our_eigenop}), which is proved in Appendix~\ref{posdef:ap},  is given by
\begin{lem}
The operator defined by (\ref{our_eigenop}) is positive-definite.
\label{posdef:lem}
\end{lem}
We shall write  $\varphi_j$, $j=0,1,2,\dots$,  for the eigenfunctions of (\ref{our_eigenop}), with corresponding eigenvalue $\lambda_j$, and normalization $\|\varphi_j\|_2=1$. Also, for the spectral problem that arises in the original problem, we shall write
\begin{equation}
\re^{\frac{1}{2}\int_0^x\Xi(s)\,\rd s}\varphi_j=x(1-x)\Psi(x)q_j.
\label{orig:eig}
\end{equation}
We shall need the asymptotic behavior of the eigenfunctions for large $\lambda_j$.

\begin{lem}
 \label{asymp:estimates}
 There exists positive constants $C_1$ and $C_2$, independent of $j$, such that
\begin{equation}
 \|\varphi_j\|_\infty\leq C_1\quad\text{and} \quad \|q_j\|_{\infty}\leq C_2\lambda_j^{3/4}.
\label{asymp:res}
\end{equation}
\end{lem}
The proof of Lemma~\ref{asymp:estimates} is given in Appendix~\ref{asymp:ap}

Finally, as in \cite{Taylor96a} for instance, we shall denote, for $s>0$, the spaces
\[
 \D_s=\left\{\left. \phi \in L^2\left([0,1],\theta\rd x\right) \right| \sum_{j=0}^\infty \widehat{\phi(j)}\lambda_j^{s/2}\varphi_j  \in L_2\left([0,1],\theta\rd x\right)\right\},\quad \widehat{\phi(j)}=(\phi,\varphi_j),
\]
with norm given by
\[
 \|\phi\|_s^2=\sum_{j=0}^\infty \widehat{\phi(j)}^2\lambda_j^s.
\]
Since Radon measures are distributions of order less or equal to zero, we have---cf. \cite{Taylor96a} with minor modifications---that:

\begin{prop}
The initial value problem defined by Equation~(\ref{our:pde:sa}) and
$w(0,x)=w^0(x)$, with $w^0\in\bmi$ has the solution
\begin{equation}
w(t,x)=\sum_{j\geq0}\widehat{w^0}(j)e^{-t\lambda_j}\varphi_j(x),\quad \widehat{w^0}(j)=(w^0,\varphi_j), 
\label{our_solution}
\end{equation}
which is unique in the class $\cicic$. 
\end{prop}
\begin{rmk}
 It can be shown that, any standard weak solution definition to (\ref{our:pde:sa}) will lead to the solution above---see for instance \cite{Evans98,Lieberman92}. Therefore, none of the conservation laws (\ref{cons:laws}) can hold, and no classical-weak solution to (\ref{our:pde}--\ref{cons:laws}) exists.
\end{rmk}

\section{Weak solution and duality formulation}

We now make precise what we mean by a weak solution to (\ref{our:pde}).
\label{formulation}

\begin{defi}
\label{our_wsd}
A weak solution to (\ref{our:pde}) will be
a function in $\lbmg$ that satisfies
\begin{align*}
&-\int_0^\infty\int_0^1p(t,x)\partial_t\phi(t,x)\d x\d t\\
&\quad=
\int_0^\infty\int_0^1p(t,x)x(1-x)\Psi(x)\left[\partial^2_{x}\phi(t,x)+\Xi(x)\partial_x\phi(t,x)\right]\d x\,\d t\\
&\qquad+\int_0^1p^0(x)\phi(0,x)\d x,
\end{align*}
where 
\[
\phi(t,x)\in \mathcal{T}=C^\infty_c\left([0,\infty)\times[0,1]\right).
\]
\end{defi}
\begin{rmk} Notice that the test functions in definition~\ref{our_wsd} are
  required to be of compact support in $[0,1]$ and not  in $(0,1)$
  as usual. Similar definitions have been given in other contexts; see
  \cite{llx01,llx06}, where they are  termed  \textsl{boundary-coupled weak solutions}.

Definition~\ref{our_wsd} can be recast in the framework   of usual
  distribution theory, by identifying a Radon measure with a compactly supported
  distribution of nonpositive order (see \cite{Tartar07}). In this case, the distribution can  
act in $C^\infty(\mathbb{R})$, but it is  entirely determined by its behavior in the support; see for instance \cite{Hormander03}.
\end{rmk}
A glance at Definition~\ref{our_wsd} shows that, on the integral on the right hand side, the test function $\phi$ is applied to the operator on the right hand side of (\ref{our:adjoint}). Thus, one could expect that any solution that satisfies (\ref{our:pde}) in the sense defined above, also satisfies the conservation laws (\ref{cons:laws}).

\begin{prop}
Let $p\in \lbmg$. If $\chi(x)$ is a stationary solution of (\ref{our:adjoint}), then the quantity 
\[
\eta(t)=\int_0^1\chi(x)p(t,x)\,\d x 
\]
is constant in time.
\end{prop}

\begin{proof}
Let $\zeta(t)\in C^{\infty}_c((0,\infty))$. Then, $\phi(t,x)=\zeta(t)\chi(x)$ is an appropriate test function. 
On substituting $\phi(t,x)$ in Definition~\ref{our_wsd}, we find that
\[
 -\int_0^{\infty}\eta(t)\zeta'(t)\rd t=0.
\]
Thus $\eta(t)$ is constant almost everywhere.
\end{proof}

\begin{rmk}
We observe that standard spectral theory shows that both the infinitesimal generators of \eqref{our:pde} and \eqref{our:adjoint} can be appropriately defined in a domain dense in $L^2((0,1))$ such that they are adjoints of each other. However, in this case, equation \eqref{our:pde} will not be the forward Kolmogorov equation associated to \eqref{our:adjoint}. On the other hand, in the sense of the pairing used in definition~\ref{our_wsd}, (\ref{our:adjoint}) with  $f(t,\cdot)\in C^{\infty}_c([0,1])$ is the adjoint of (\ref{our:pde}) with $p(t,\cdot)\in\bmg$. Thus, we recover the usual interpretation of the conservation laws given by the kernel of the adjoint.
\end{rmk}

\section{Existence and uniqueness}
\label{the:theorem}

In what follows, it will be convenient to decompose a compact distribution, or a Radon measure, as the sum of a distribution without  singular support at the endpoints, and two distributions  singularly supported at the endpoints. We shall write $\E$ to denote the space of compactly supported distributions in $\rb$.
\begin{lem}
Let $\nu\in\E$, with $\supp(\nu)=[0,1]$. Then, the setwise decomposition
\[
 [0,1]=\{0\}\cup(0,1)\cup\{1\},
\]
yields a decomposition of $\nu$, namely
\[
 \nu=\nu_0+\mu+\nu_1,
\]
where $\nu_i$ is a compact distribution supported at $x=i$, and  $\ssupp(\mu)\subset(0,1)$. Moreover, if $\nu$ is a Radon measure, then $\mu\in\bmgi$, and $\nu_i=c_i\delta_i$ , with $c_i\in\rb$,  are singular measures with support at $x=i$.
\end{lem}
\begin{proof}
Let $\zeta_i^{\eps}$, $i=0,1,2$ be a partition of unity in $[0,1]$, subordinated to the open cover $\left\{[0,2\eps),  (1-2\eps,1], (\eps,1-\eps)\right\}$. Let $\phi\in C^{\infty}_c([0,1])$. Define $\nu_i$, $i=0,1$ and $\mu$ by
\begin{align*}
 \int_0^1\nu_i\phi(x)\,\rd x & :=\lim_{\eps\to0}\int_0^1\zeta_i^{\eps}\nu\phi(x)\,\rd x,\quad i=0,1\ ,\\
\int_0^1\mu\phi(x)\,\rd x & :=\lim_{\eps\to0}\int_0^1\zeta_2^{\eps}\nu\phi(x)\,\rd x.\\
\end{align*}
Then clearly $\nu= \nu_0+\mu+\nu_1$.  Also, it is readily seen that  $\ssupp(\mu)\subset (0,1)$. Since $\zeta_0^{\eps}(x)=1$ and ${\zeta_0^{\eps}}^{(n)}(x)=0$, $n\geq1$, for $x\in[0,\eps)$, we find that $\nu_0$ is supported at $x=0$, with a similar argument holding for $\nu_1$.  Moreover, since a Radon measure is inner regular, the restriction of $\nu$ to $(0,1)$ yields a Radon measure in $(0,1)$. Finally, a Radon measure supported in a singleton must be an atomic measure.
\end{proof}

For the initial condition, we shall write
\[
 p^0=a^0\delta_0+ q^0 + b^0\delta_1,
\]
to denote the corresponding decomposition.
Also, in order to show the existence of a solution to (\ref{our:pde}) in the sense of definition~\ref{our_wsd}, we shall temporarily consider $p\in\len$, with support in $[0,1]$.
We shall write
\[
 p=p_0+q+p_1,
\]
for the decomposition of $p$.

We now show that $q$ must be, as a matter of fact, much more regular. 
\begin{prop}[Interior regularity]
\label{intreg}
 Assume that $q^0\in\bmi$. If a solution to (\ref{our:pde}) exists, then $q(t,x)$ 
must be the unique  classical solution in the sense of section~\ref{prelim}, 
with $q(0,x)=q^0(x)$, i.e.,
\begin{equation}
 q(t,x)=\sum_{j=0}^{\infty}\widehat{q^0}(j)q_j\re^{-\lambda_jt},\label{qform}
\end{equation}
with  $q_j$ given by \eqref{orig:eig}, and 
$\widehat{q^0}(j)$ is $j$-th Fourier coefficient of $q^0$.
In particular, 
\[
q\in\cicic\ .
\]
\end{prop}
\begin{proof}
 Let $\phi \in C^{\infty}_c\left([0,\infty)\times(0,1)\right)$. Applying to definition~\ref{our_wsd}, we find
\begin{align*}
&-\int_0^{\infty}\int_0^1q(t,x)\partial_t\phi(t,x)\,\rd x\,\rd t\\
&\quad =
\int_0^{\infty}\int_0^1q(t,x)x(1-x)\Psi(x)
\left[\partial_x^2\phi(t,x)+\Xi(x)\partial_x\phi(t,x)\right]\,\rd x\,\rd t \\
&\qquad+ \int_0^{\infty} q^0(x)\phi(0,x)\,\rd x.
\end{align*}
The result now follows by taking testing functions of the form
\[
\phi(t,x)=\zeta(t)\re^{-\frac{1}{2}\int_0^x\Xi(s)\,\rd s}\tilde{\phi}(x),\quad 
\zeta\in C_c([0,\infty))\quad\text{and}\quad \tilde{\phi}\in C_c((0,1)),
\]
and then we use a standard Galerkin approximation procedure.
\end{proof}
Before we proceed, we observe that, since  $p_0$ and $p_1$  are distributions supported on a singleton, we must have, for some integers $M$ and $M'$, that 
\begin{equation}
p(t,x)=
\sum_{k=0}^Ma_k(t)\delta_0^{(k)}+\sum_{k=0}^{M'}b_k(t)\delta_1^{(k)} + q(t,x),\label{weak_form} 
\end{equation}
where $\delta_{x_0}^{(k)}$ denotes the $k$-th distributional derivative of the singleton supported measure.

\begin{thm}[Existence and uniqueness]
 The unique solution of (\ref{our:pde}) in the sense of definition~\ref{our_wsd}, with initial condition $p^0\in\bm$ is given by
\[
 p(t,x)=q(t,x)+a(t)\delta_0+b(t)\delta_1,
\]
with $q(t,x)$ given by (\ref{qform}). Moreover, we have
\[
a(t)=\Psi(0)\int_0^tq(s,0)\d s + a^0\quad\text{and}\quad
b(t)=\Psi(1)\int_0^tq(s,1)\d s + b^0.
\]
\end{thm}
\begin{proof}

First, we define
\[
 \bt=\left\{\phi\in C^{\infty}_c\left((0,\infty)\times[0,1]\right)\right\}.
\]
For $l>0$, we also define,
\[
 \bti{l}{0}=\left\{\phi\in C_c([0,\infty)\times[0,1))|\partial_x^n\phi(t,0)=0,0\leq n<l\right\},
\]
with a similar definition for $\bti{l}{1}$. Notice that, for $r>s$, 
$\bti{r}{0}\subset\bti{s}{0}$.

On substituting (\ref{weak_form}) in definition~\ref{our_wsd}, with $\phi\in\bt$, using that $q$ is smooth for $t>0$ and integrating by parts we obtain that
\begin{align*}
-\int_0^{\infty}\left[\sum_{k=0}^Ma_k(t)\partial_t\partial_x^k\phi(t,0) +
\sum_{k=0}^{M'}b_k(t)\partial_t\partial_x^k\phi(t,1)\right]\,\rd t=
\int_0^{\infty}\left[q(t,0)\phi(t,0)+q(t,1)\phi(t,1)\right]\,\rd t +\\
+\int_0^{\infty}\sum_{k=0}^Ma_k(t)\partial_x^k\left.(x(1-x)\Psi(x)\partial_x^2\phi(t,x))\right|_{x=0}\,\rd t
+ 
\int_0^{\infty}\sum_{k=0}^Ma_k(t)\partial_x^k\left.(x(1-x)\Pi(x)\partial_x\phi(t,x))\right|_{x=0}\,\rd t+\\
+\int_0^{\infty}\sum_{k=0}^{M'}b_k(t)\partial_x^k\left.(x(1-x)\Psi(x)\partial_x^2\phi(t,x))\right|_{x=1}\, \rd t
+\int_0^{\infty}\sum_{k=0}^{M'}b_k(t)\partial_x^k\left.(x(1-x)\Pi(x)\partial_x\phi(t,x))\right|_{x=1}\,\rd t.
\end{align*}
 Restricting somewhat further, for $\phi\in\bti{M+1}{0}$, we find that
\[
 0=\int_0^\infty
a_M(t)\partial_x[x(1-x)\Psi(x)]_{x=0}\partial_x^{M+1}\phi(t,0)\,\rd t.
\]
Thus $a_M(t)=0$, and the sum can be only up to $M-1$. Repeating the argument inductively yields $M=0$. An analogous argument yields $M'=0$. Thus only $a_0(t)$ and $b_0(t)$ can be nonzero. We now drop the subscripts and determine their values.
Applying definition~\ref{our_wsd} to $\phi\in\bt$, such that $\phi(t,1)=0$, we find that
\[
-\int_0^{\infty}a(t)\partial_t\phi(t,0)\,\rd t = \int_0^{\infty}q(t,0)\Psi(0)\phi(t,0)\,\rd t .
\]

Integrating by parts the corresponding relation for $a(t)$, we obtain 
\[
\int_0^\infty a(t)\partial_t\phi(t,0)\,\rd t=
\int_0^\infty \left(\Psi(0) \int_0^tq(s,0)\,\rd s + a^0\right) \partial_t\phi(t,0)\,\rd t .
\]
Hence
\[
a(t)-\Psi(0)\int_0^tq(s,0)\,\rd s - a^0=\mathrm{const},
\]
everywhere, in as much as the integral is continuous. Since $a(0)=a^0$, the identity follows. A similar calculation also shows that
\[
b(t)=\Psi(1)\int_0^tq(s,1)\,\rd s +b^0.
\]
Uniqueness follows from proposition~\ref{intreg} and from the expressions for $a(t)$ and $b(t)$.
Finally, notice that, since $q(t,x)\geq0$, we have that both $a$ and $b$ are increasing.
\end{proof}

\section{Large time behavior}
\label{largetime}

We now present some  results  for the behavior of the solution in the large time limit.

Let us define
\begin{align*}
 b^{\infty}&:=\int_0^1\psi(x)p^0(x)\,\rd x = b^0 + \int_0^1 \psi(x)q^0(x)\,\rd x,\\
a^{\infty}&:=\int_0^1p^0(x)\,\rd x -b^{\infty} = a^0+\int_0^1(1-\psi(x))q^0(x)\,\rd x.
\end{align*}
Using the conservation laws \eqref{cons:laws}, we also have
\begin{align*}
 b^{\infty}&:=\int_0^1\psi(x)p(t,x)\,\rd x = b(t)+ \int_0^1 \psi(x)q(t,x)\,\rd x,\\
a^{\infty}&:=\int_0^1p(t,x)\,\rd x -b^{\infty} = a(t)+\int_0^1(1-\psi(x))q(t,x)\,\rd x.
\end{align*}
Since $0\leq\psi\leq1$ and $q(t,\cdot)\geq0$, we have that both $a^{\infty}-a(t)$ and $b^{\infty}-b(t)$ are nonnegative. From the representation given by \eqref{our_solution}, we have that $\lim_{t\to\infty}\|q(t,x)\|_{\infty}=0$. Hence, $\lim_{t\to\infty}a(t)=a^{\infty}$ and $\lim_{t\to\infty}b(t)=b^{\infty}$.

Moreover, since $q(t,x)\geq0$, we have
\[
 a^{\infty}-a(t)+b^{\infty}-b(t)=\int_0^1q(t,x)\,\rd x=\|q(t,\cdot)\|_1, 
\]
which also yields the inequalities
\[
 a^{\infty}-a(t) \leq \|q(t,\cdot)\|_1\quad\text{and}\quad b^{\infty}-b(t)\leq \|q(t,\cdot)\|_1.
\]
The behavior of the $L^1$ norm of $q$ is given by the following result:

\begin{prop}
 \label{estimates}
Let $p$ be the solution to (\ref{our:pde}) with an initial condition with $q^0\in\bmi$ and let $\lambda_0$ be the smallest eigenvalue of (\ref{our_eigenop}). Then we have that
\[
 \lim_{t\to\infty}\re^{\lambda_0t}\|q(t,\cdot)\|_1=C_{\infty}.
\]
In addition, if we assume that
\[
 w^0=x(1-x)\Psi(x)\re^{-\frac{1}{2}\int_0^x\Xi(s)\,\rd s}q^0\in \bmi\cap \mathcal{D}_s,\quad s>0,
\]
then there exists $C_{0,s}>0$ such that
\[
 ||q(t,\cdot)||_1\leq C_{0,s}\|w^0\|_s\re^{-\lambda_0t}.
\]
In particular, the same limit property and bounds apply to $a^{\infty}-a(t)$ and $b^{\infty}-b(t)$.
\end{prop}
\begin{proof}
For the first part, recall that
\[
 q(t,x)=\sum_{j=0}^{\infty}\widehat{w^0}(j)\re^{-\lambda_jt}q_j(x).
\]
Let us write
\[
 Q_j=\int_0^1q_j(x)\,\rd x .
\]
We observe that
\[
 Q_j=(\re^{-\frac{1}{2}\int_0^x\Xi(s)\,\rd s},\varphi_j).
\]
However, since $\re^{\frac{1}{2}\int_0^x\Xi(s)\,\rd s}\not\in L^2([0,1],\theta\rd x)$, we do not have an immediate bound for $|Q_j|$.

On the other hand, we observe that $q_j$ satisfies
\[
 -\lambda_jq_j(x)=\partial_x^2\left[x(1-x)\Psi(x)q_j(x)\right]-\partial_x\left[x(1-x)\Pi(x)q_j(x)\right],
\]
which integrated yields
\[
 Q_j = \frac{\Psi(0)q_j(0)+\Psi(1)q_j(1)}{\lambda_j}.
\]
For large $j$, \eqref{asymp:res} guarantees that we must then have 
\[
 |Q_j|\leq C_2\lambda_j^{-1/4}.
\]
Thus, for $t>0$, we have 
\[
 \re^{\lambda_0t}\|q(t,\cdot)\|_1=Q_0\widehat{w^0}(0)+\sum_{j=1}^\infty Q_j\widehat{w^0}(j)\re^{-(\lambda_j-\lambda_0)t},
\]
and the result follows with $C_{\infty}=Q_0\widehat{w^0}(0)$.

For the second part, let us define the auxiliary functions:
\[
 \alpha_s(x)=\sum_{j=0}^\infty \widehat{w^0(j)}\lambda_j^{s/2}\varphi_j(x)
\quad\text{and}\quad
\beta_s(t,x)=\sum_{j=0}^\infty Q_j\lambda_j^{-s/2}\varphi_j(x)\re^{-\lambda_jt}.
\]
Then, we have that
\begin{align*}
 \|q(t,\cdot)\|_1&=(\alpha_s,\beta_s(t,\cdot))\leq \|\alpha_s\|_2\|\beta_s(t,\cdot)\|_2=\\
&=\|w^0\|_s\re^{-\lambda_0t}\|\re^{\lambda_0t}\beta_s(t,\cdot)\|_2\leq C_{0,s}
\|w^0\|_s\re^{-\lambda_0t},
\end{align*}
with $C_{0,s}=\|\beta_s(0,\cdot)\|_2$.

\end{proof}

\begin{thm}[Exponential convergence]
Let $\rho$ denote the Radon metric, and let
\[
 p^{\infty}=a^{\infty}\delta_0+b^{\infty}\delta_1.
\]
Under the same hypothesis of proposition~\ref{estimates}, we have that 
\begin{equation}
 \lim_{t\to\infty}\re^{\lambda_0t}\rho(p,p^{\infty})\leq 2C_{\infty}.\label{radonconvergence}
\end{equation}
With the additional hypothesis, we have that
\begin{equation}
 \rho(p,p^{\infty})\leq 2C_{0,s}\|w^0\|_s\re^{-\lambda_0t}.\label{radonbound}
\end{equation}

In particular, (\ref{radonconvergence}) implies convergence in the Wasserstein metric.
\end{thm}
\begin{proof}
 Recall that 
\[
 \rho(\nu,\mu)=\sup\left\{\left.\int_0^1f(x)\rd(\nu-\mu)\,\right| \, f\in C([0,1];[-1,1])\right\}.
\]
But, for such $f$ we have that, when $t>0$.
\begin{align*}
 \left|\int_0^1f(x)\rd(p^{\infty}-p(t,x))\right|& \leq  \int_0^1|\rd(p^{\infty}-p(t,x))|\\
&\leq \int_0^1(a^{\infty}-a(t))\delta_0\rd x+\int_0^1(b^{\infty}-b(t))\delta_1\rd x + \int_0^1|q(t,x)|\rd x\\
&=a^{\infty}-a(t) + b^{\infty}-b(t) + \|q(t,\cdot)\|_1\\
&= 2\|q(t,\cdot)\|_1.
\end{align*}
Now, both \eqref{radonconvergence} and \eqref{radonbound} follows from Proposition~4.
\end{proof}

\begin{rmk}
 In many applications, the slowest decaying mode $\varphi_0$ is taken to be a quasi-stationary distribution for the diffusion process. The constant $C_{\infty}$ is then the  total probability mass of such a distribution. 
\end{rmk}

\appendix

\section{Postponed proofs}

\subsection{Proof of positive-definiteness of \eqref{our_eigenop}}
\label{posdef:ap}
\begin{proof}[Proof of Lemma~\ref{posdef:lem}]
Let
\begin{equation}
 v=\re^{-\frac{1}{2}\int_0^x\Xi(s)\,\rd s}\varphi.
\label{our_change_adj}
\end{equation}
Then \eqref{our_eigenop} becomes
\begin{equation}
  \begin{array}{c}
  -v''-\Xi v'=\lambda \theta(x)v, \qquad  v(0)=v(1)=0.
 \end{array}
\label{our_eigenop_adj}
\end{equation}
When $\lambda=0$, then \eqref{our_eigenop_adj} becomes the stationary version of \eqref{our:adjoint2}, with Dirichlet boundary condition. Its  general solution, $\bar{v}$, is given by $\bar{v}=c_1 + c_2\psi$, which does not satisfy the required boundary conditions. Thus, zero cannot be an eigenvalue of \eqref{our_eigenop_adj}. 
Moreover, since the transformation \eqref{our_change_adj} preserves the oscillation properties of the eigenfunctions, we have that the eigenfunction $v_0$, corresponding to the smallest eigenvalue $\lambda_0$, will not have any zeros inside $(0,1)$. Let us assume, without loss of generality, that $v_0>0$ in $(0,1)$. It must have a point of maximum $x=x^*\in(0,1)$, where $v_0'(x^*)=0$. Hence we must have
\[
 -\lambda_0\frac{v_0(x^*)}{x^*(1-x^*)\Psi(x^*)}=v_0''(x^*).
\]
Note that $v_0''(x^*)\not=0$, otherwise we would have $\lambda_0=0$. Since it is a maximum, we must have $v_0''(x^*)<0$. Since, $v_0(x^*)>0$, we  have $\lambda_0>0$.
\end{proof}

\subsection{Proof of the asymptotic estimates}
\label{asymp:ap}

\begin{proof}[Proof of Lemma~\ref{asymp:estimates}]

For the proof, we drawn on results by \cite[chapter 12]{Olver74} that are summarized as follows

\begin{thm}
\label{asymp:stuff}
Let
\[
 \frac{1}{\zeta}\left(\frac{\d\zeta}{\d x}\right)^2=-4\theta(x),\quad \zeta(0)=0.
\]
Also let
\[
 \hvp_j(\zeta)=A_{0,j}\left(\frac{\d\zeta}{\d x}\right)^{-1/2}|\zeta|^{1/2}J_1\left(\lambda_j^{1/2}|\zeta|^{1/2}\right),
\]
where $J_1$ is the standard Bessel function of order one, and $A_{0,j}$ is choosen such that $\|\hvp_j\|_2=1$. For large $j$, we have that
\[
 \|\varphi_j -\hvp_j\|_{\infty}\leq 
K_1F_1(\lambda_j^{1/2}|\zeta|^{1/2})\exp\left(\frac{K_2}{\lambda_j^{1/2}}F_2(\zeta)\right)\frac{F_2(\zeta)}{\lambda_j^{1/2}},
\]
where $K_1$, $K_2$ are positive constants and  $F_1$, $F_2$ are positive and bounded continuous functions.
\end{thm}
With this result,  we can now prove the asymptotic behavior for $\varphi_j$ and $q_j$

Let
\[
 A_{0,j}^{-2}=\int_0^1\left(\frac{\d\zeta}{\d x}\right)^{-1}|\zeta|J_1^2\left(\lambda_j^{1/2}|\zeta|^{1/2}\right) \,\rd x
\]
and $J_1$ is the standard Bessel function of order one.
Let $z=\lambda_j^{1/2}|\zeta|^{1/2}$. Then we find that
\[
 \int_0^1\left(\frac{\d\zeta}{\d x}\right)^{-1}|\zeta|J_1^2\left(\lambda_j^{1/2}|\zeta|^{1/2}\right) \theta\,\rd x=
\frac{1}{\lambda_j}\int_0^{z_1}zJ_1^2(z)\, \rd z,
\]
where $z^1=\lambda_j^{1/2}|\zeta(1)|$. For large $j$, we have from the asymptotic behavior of $J_1$ at infinity that 
\[
 A_{0,j}=C\lambda_j^{1/4} + \O(1).
\]
Also, since $uJ_1(u)\leq C\sqrt(u)$, for large $u$, we have 
\[
 \left\||\zeta|^{1/2}J_1\left(\lambda_j^{1/2}|\zeta|^{1/2}\right)\right\|_{\infty}=\frac{1}{\lambda_j^{1/2}}\|zJ_1(z)\|_{\infty}\leq
C\lambda_j^{-1/4}.
\]
Combining these two results, with Theorem~\ref{asymp:stuff}, we have the first result in \eqref{asymp:res}. 

For the second one, we observe that 
\[
 \|\theta\hvp_j\|_{\infty}=A_{0,j}\lambda_j^{1/2}\left\|\left(\frac{\d\zeta}{\d x}\right)^{3/2}\frac{J_1(z)}{z}\right\|_{\infty}\leq CA_{0,j}\lambda_j^{1/2}.
\]
Combining with the estimate for $A_{0,j}$, we have the result.
\end{proof}
%


\end{document}